\documentclass{amsart}

\usepackage{amsmath,amstext,amsfonts,amssymb}

\begin{document}
\title{Toeplitz-Composition $C^*$-Algebras}

\author{Thomas L. Kriete}
\address{Department of Mathematics,
University of Virginia, Charlottesville, VA 22904}
\email{tlk8q@virginia.edu}

\author{Barbara D. MacCluer}
\address{Department of Mathematics,
University of Virginia, Charlottesville, VA   22904}
\email{bdm3f@virginia.edu}

\author{Jennifer L. Moorhouse}
\address{Department of Mathematics,
Colgate University, Hamilton, NY   13356}
\email{jmoorhouse@mail.colgate.edu}


\date{July 14, 2005}

\newtheorem{theorem}{Theorem}
\newtheorem{prop}{Proposition}
\newtheorem{lemma}{Lemma}
\newtheorem{cor}{Corollary}

\newcommand{\p}{\varphi}
\newcommand{\cp}{C_{\varphi}}
\newcommand{\D}{{\mathbb D}}
\newcommand{\B}{{\mathcal B}}
\newcommand{\com}{{[C_{\varphi}^*,C_{\varphi}]}}

 \newcommand{\ebox}{\hspace*{1em}\linebreak[0]\hspace*{1em}\hfill\rule[-1ex]{
.1in}{.1in}\\}
 \newcommand{\pf}{\noindent \bf Proof.\rm\ \ }
 
\maketitle

\begin{center}
{\it In memory of Marvin Rosenblum}
\end{center}

\begin{abstract}
Let $\zeta$ and $\eta$ be distinct points on the unit circle and suppose
that $\varphi$ is a linear-fractional self-map of the unit disk $\D$, not
an automorphism, with $\varphi(\zeta)=\eta$.  We describe the $C^*$-algebra
generated by the associated composition operator $\cp$ and the 
shift operator,
acting on the Hardy
space on $\D$.
\end{abstract}

\maketitle

\section{INTRODUCTION}

Any analytic self-map $\varphi$ of the unit disk $\D$ induces a 
bounded composition
operator $\cp:f\rightarrow f\circ\varphi$ on the Hardy space $H^2$.  
The linear-fractional self-maps of $\D$
form a rich class of examples, and many properties of composition
operators are profitably studied in the context of these maps
(e.g. cyclicity, spectral properties, subnormality; 
see \cite{Co}, \cite{CMbook}, \cite{Sbook}).  The space $H^2$ also 
supports the Toeplitz operators
$T_w$.  Here, $w$ is a bounded measurable function on the unit 
circle $\partial \D$,
and $T_w$ acts on $H^2$ by $T_wf=P(wf)$, where $P$ is the orthogonal 
projection
of $L^2$ (the Lebesgue space associated with normalized arc-length 
measure on $\partial\D$)
onto $H^2$. Taking $w$ to be the independent variable $z$, one obtains 
the shift 
operator $T_z$ on $H^2$.  
A theorem of I. Gohberg and 
I. Fel'dman (\cite{G}, \cite{GF}) and
L. Coburn (\cite{Cob1}, \cite{Cob2}) asserts  
that $C^*(T_z)$, the
unital $C^*$- algebra generated by $T_z$, contains the 
ideal $\mathcal{K}$ of compact operators,
as well as all Toeplitz operators $T_w$ with continuous 
symbol $w$.  Moreover, the map sending
$w$ to the coset of $T_w$ is a $*$-isomorphism of 
$C(\partial\D)$, the algebra
of continuous functions on $\partial \D$, onto the 
quotient algebra $C^*(T_z)/\mathcal{K}$.
In this article our goal
is to replace $C^*(T_z)$ by $C^*(T_z,\cp)$, the unital 
$C^*$-algebra generated by $T_z$ and
$\cp$, for certain linear-fractional $\varphi$.  

Section 2 presents a characterization of those analytic self-maps $\varphi$ of $\D$
with $|\varphi(e^{i\theta})|<1$ a.e. on $\partial \D$ for which $\cp$ commutes with
$T_z$ or $T_z^*$ modulo $\mathcal{K}$.  
In Section 3 we show that for any linear-fractional self-map $\varphi$ 
of the disk
which is not an automorphism, there is an associated linear-fractional 
map $\sigma$
(the ``Krein adjoint" of $\varphi$) and a scalar $s$ so that 
$\cp^*=sC_{\sigma}+K$ for some compact operator $K$.  Our setting here
is primarily that of $H^2$, although this result is easily extended to
the Bergman space.
This theorem plays a key role in the
work in Section 4, where we study $C^*(T_z,\cp)$. Recent work of 
M. Jury \cite{J} treats the case where
$\varphi$ is an automorphism (and indeed ranges over a discrete group $\Gamma$
of automorphisms), showing that the $C^*$-algebra generated by 
$\{C_{\varphi}:\varphi\in \Gamma\}$
contains $T_z$, and exhibiting the quotient of this algebra by $\mathcal{K}$
as the discrete crossed product $C(\partial\D)\times\Gamma$.  In the
present article we suppose $\varphi$ is not an automorphism
but does satisfy $\|\varphi\|_{\infty}=1$.  In the case that $\varphi$ 
is a parabolic 
non-automorphism (see Section 2
for a discussion of this terminology; such maps have a fixed point 
on $\partial \D$),
the work of P. Bourdon, D. Levi, S. Narayan and J. Shapiro 
in \cite{BLSS} shows that $\cp^*\cp-\cp\cp^*$ is compact.  Such a $\cp$ 
also commutes with 
$T_z$ and $T_z^*$ modulo
$\mathcal{K}$, so that $C^*(T_z,\cp)/\mathcal{K}$ is commutative, 
hence describable
by Gelfand theory.  Here we suppose that $\varphi$ is neither an 
automorphism
nor a parabolic non-automorphism, but that there exist distinct points
$\zeta,\eta$ in $\partial\D$ with $\varphi(\zeta)=\eta$.  In this case
$C^*(T_z,\cp)/\mathcal{K}$ is not commutative, but we will see that
it is tractable.  As an application, in Section 4.6 we concretely 
determine the essential spectrum of
any element of $C^*(T_z,\cp)$.  Our main tool is the localization 
theorem of R. G. Douglas \cite{Do}.

We thank Paul Bourdon for several helpful comments.

\section{COMPOSITION OPERATORS ESSENTIALLY COMMUTING WITH $T_z$ OR $T_z^*$}

The commutator $AB-BA$ of two bounded operators $A$ and $B$ on a Hilbert
space $\mathcal{H}$ is denoted $[A,B]$. 
An operator is said to be essentially normal if its self-commutator 
$[A^*,A]$ is compact.
In the course of their work on essentially normal linear-fractional 
composition
operators, Bourdon, Levi, Narayan and Shapiro \cite{BLSS} show that 
if $\varphi$ is 
a linear
fractional non-automorphism mapping $\D$ into $\D$
and fixing a point of $\partial \D$, then $[T_z^*,\cp]$ is compact, where
$T_z$ is the shift on $H^2$.
Here we will give a generalization which is perhaps of independent interest.

For $\alpha$ a complex number of modulus $1$, and $\varphi$ an analytic
self-map of $\D$, the real part of $(\alpha+\varphi)/
(\alpha-\varphi)$ is a positive harmonic function on $\D$.  Necessarily then 
this function
is the Poisson integral of a finite positive Borel measure $\mu_{\alpha}$
on $\partial\D$; $\mu_{\alpha},|\alpha|=1$ are the {\it Clark measures}
for $\varphi$.  We write $E(\varphi)$ for the closure in $\partial \D$ of the union
of the closed supports of the singular parts $\mu_{\alpha}^s$ of the Clark measures
as $\alpha$ ranges over the unit circle.  For a linear-fractional non-automorphism
$\varphi$ which 
sends $\zeta\in\partial\D$ to $\eta\in\partial\D$, one has $\mu_{\alpha}^s=0$
when $\alpha\neq \eta$ and $\mu_{\eta}^s=|\varphi'(\zeta)|^{-1}\delta_{\zeta}$,
where $\delta_{\zeta}$ is the unit point mass at $\zeta$.
We will use the following result, proved in \cite{KM}.  Here $M_w$ denotes
the operator on $L^2=L^2(\partial\D)$
of multiplication by the bounded measurable function $w$.  

\begin{theorem}\label{fromKM}\cite{KM}
Let $\varphi$ be an analytic self-map of $\D$ such that $|\varphi(e^{i\theta})|<1$
a.e. with respect to Lebesgue measure on $\partial \D$, and 
suppose that $w$ is a 
bounded measurable function on $\partial \D$ which is continuous at each
point of $E(\varphi)$.  The weighted composition operator 
$M_w\cp:H^2\rightarrow L^2$
is compact if and only if $w\equiv 0$ on $E(\varphi)$.
\end{theorem}

It will be convenient to recast Theorem~\ref{fromKM} in terms of Toeplitz operators.

\begin{cor}\label{thm1cor}
Suppose that $\varphi$ and $w$ satisfy the hypotheses in the first sentence of
Theorem~\ref{fromKM}.  Then $T_w\cp:H^2\rightarrow H^2$ is compact if 
and only if $w\equiv 0
$ on $E(\varphi)$.  
\end{cor}
\begin{proof}
It is enough to show that $M_w\cp$ is compact when $T_w\cp$ is compact.
Note that 
$$M_w\cp=T_w\cp+H_w\cp$$
where $H_w:H^2\rightarrow(H^2)^{\perp}$ is the Hankel operator defined by
$H_w=(I-P)M_w|_{H^2}$.  We need only check that $H_w\cp$ is compact.
Let $\tilde{w}$ be a continuous function on $\partial \D$ agreeing with
$w$ on $E(\varphi)$.  We have
$$H_w\cp=(I-P)M_{(w-\tilde{w})}\cp+H_{\tilde{w}}\cp.$$
Since $\tilde{w}$ is continuous, $H_{\tilde{w}}$ is compact by Hartman's theorem
\cite{Ha}.  On the other hand, $M_{(w-\tilde{w})}\cp$ is compact by 
Theorem~\ref{fromKM},
and we are done.
\end{proof}

The next result gives the above-mentioned generalization.  

\begin{theorem}\label{commutes}
Let $\varphi$ be an analytic self-map of $\D$ such that 
$|\varphi(e^{i\theta})|<1$
a.e. with respect to Lebesgue measure.  Suppose that $\varphi$ agrees 
almost everywhere on $\partial \D$ with a bounded measurable function 
$\widehat{\varphi}$ 
which is continuous
at each point of $E(\varphi)$.  Then the following are equivalent:
\begin{itemize}
\item[(i)] $[T_z,\cp]\in\mathcal{K}$.
\item[(ii)] $[T_z^*,\cp]\in\mathcal{K}$.
\item[(iii)] For each $\zeta$ in $E(\varphi)$, 
$\widehat{\varphi}(\zeta)=\zeta$.
\end{itemize}
When these conditions hold, $[T_w,\cp]\in\mathcal{K}$ for every 
$w$ in $C(\partial\D)$.
\end{theorem}
\begin{proof}
We use the following identity from \cite{BLSS}:
\begin{equation}\label{backshift}
[T_z^*,\cp]=T_{(\overline{z}\varphi-1)}\cp T_z^*.
\end{equation}
Since $T_z^*$, the backward shift, is a partial isometry with range
$H^2$, the operator on the right-hand side of Equation~(\ref{backshift})
is compact exactly when $T_{(\overline{z}\varphi-1)}\cp$ is compact.
This operator clearly coincides with 
$T_{(\overline{z}\widehat{\varphi}-1)}\cp$. Corollary~\ref{thm1cor}
gives the equivalence of (ii) and (iii). 
For the equivalence of (i) and (iii) we easily
check that 
$$[T_z,\cp]=T_{(z-\varphi)}\cp=T_{(z-\widehat{\varphi})}\cp$$
and again apply Corollary~\ref{thm1cor}, with $w=z-\widehat{\varphi}$. 
The statement about $[T_w,\cp]$ is immediate.
\end{proof}

\section{THE ADJOINT OF $\cp$}

In this section we develop some properties of linear-fractional composition
operators and their adjoints.
To any linear-fractional map 
\begin{equation}\label{varphi}
\varphi(z)=(az+b)/(cz+d)
\end{equation}
we associate another linear-fractional
map $\sigma_{\varphi}$ defined as 
\begin{equation}\label{defofsigma}
\sigma_{\varphi}(z)=(\overline{a}z-\overline{c})/(-\overline{b}z+\overline{d}).
\end{equation} 
The map $\sigma_{\varphi}$ is sometimes referred to as 
the ``Krein adjoint" of $\varphi$; for
an explanation of this terminology, see \cite{CM}.
When no confusion can result, we write $\sigma$ for $\sigma_{\varphi}$.
When $\varphi$ is a self-map of the disk, $\sigma$ will be also, 
and if $\varphi(\zeta)=\eta$
for $\zeta,\eta\in\partial\D$, then $\sigma(\eta)=\zeta$; see \cite{Co}.
Carl Cowen \cite{Co} has shown that the adjoint of any 
linear-fractional $\cp$, acting on $H^2$, is given by
\begin{equation}\label{cowenadjoint}
\cp^*=T_gC_{\sigma}T_h^*
\end{equation}
where $g(z)=(-\overline{b}z+\overline{d})^{-1}$, $h(z)=cz+d$, 
and $T_g, T_h$ are
the analytic Toeplitz operators of multiplication by the $H^{\infty}$ 
functions $g$ and $h$.

Our first result uses Equation~(\ref{cowenadjoint}) to show that
when $\|\varphi\|_{\infty}=1$ but $\varphi$ is not an automorphism, 
the adjoint of $\cp$, modulo 
the ideal ${\mathcal K}$ of compact operators, is a
scalar multiple of $C_{\sigma}$.

\begin{theorem}\label{adjointcp}
Suppose that $\varphi$ given by Equation~(\ref{varphi}) is a 
linear-fractional self-map of $\D$, not an automorphism,
which satisfies $\varphi(\zeta)=\eta$ for some $\zeta,\eta\in\partial \D$.
Let 
$s=(\overline{c}\overline{\zeta}+\overline{d})/(-\overline{b}\eta+\overline{d})$.  
Then
there exists a
compact operator $K$ on $H^2$ so that $\cp^*=sC_{\sigma}+K$, where
$\sigma$ is as given by Equation~(\ref{defofsigma}). 
\end{theorem}
\begin{proof}
We first consider the case where $\zeta=\eta$, so that $\zeta$ is a fixed
point of $\varphi$. Let $\sigma, h$ and $g$ be associated to
$\varphi$ as in Equations~(\ref{defofsigma}) and (\ref{cowenadjoint}), and note
that $\sigma$ fixes $\zeta$ also.  It is immediate
that
$[C_{\sigma},T_h^*]=\overline{c}[C_{\sigma},T_z^*]$.  Invoking 
Theorem~\ref{commutes},
it follows that $C_{\sigma}T_h^*=T_h^*C_{\sigma}+K_1$ for some
compact operator $K_1$.  From
Equation~(\ref{cowenadjoint}) we then have
\begin{eqnarray*}
\cp^*&=&T_gC_{\sigma}T_h^*\\
&\equiv& T_gT_h^*C_{\sigma}\  (\mbox{mod }\mathcal{K})\\
&\equiv&T_{\overline{h}g}C_{\sigma}\  (\mbox{mod }{\mathcal K})
\end{eqnarray*}
where the last line is justified by Proposition 7.22 in \cite{Do}.
Since $E(\sigma)=\{\eta\}=\{\zeta\}$, we may now apply Corollary~\ref{thm1cor}
with $w=\overline{h}g-\overline{h(\zeta)}g(\zeta)$ to see that
$$T_{\overline{h}g}C_{\sigma}-\overline{h(\zeta)}g(\zeta)C_{\sigma}=
T_{(\overline{h}g-
\overline{h(\zeta)}g(\zeta))}C_{\sigma}\in\mathcal{K},$$
which gives the desired conclusion.

In the case that $\zeta\neq \eta$ we consider the map 
$\psi(z)=\zeta\overline{\eta}\varphi(z)$
which fixes $\zeta$.  Since $C_{\varphi}^*=C_UC_{\psi}^*$ 
where $U(z)=\zeta\overline{\eta}z$,
the first part of the argument shows that $C_{\varphi}^*=C_UC_{\psi}^*\equiv
C_U(\overline{h}_{\psi}(\zeta)g_{\psi}(\zeta)C_{\sigma_{\psi}})
\ (\mbox{mod }{\mathcal K})$.  Since $\sigma_{\psi}\circ U=\sigma_{\varphi}$ 
and $\overline{h}_{\psi}(\zeta)
g_{\psi}(\zeta)=
(\overline{c}\overline{\zeta}+\overline{d})/(-\overline{b}\eta+\overline{d})$,
the conclusion follows.
\end{proof}

\noindent {\bf Remark 1.} An analogue of Theorem~\ref{adjointcp} holds 
in the Bergman
space $A^2$ of analytic functions in $L^2(\D,dA)$, where $dA$ is normalized area
measure on $\D$.  If $\varphi$ given by Equation~(\ref{varphi}) is a self-map
of $\D$, then on $A^2$ we have $\cp^*=T_gC_{\sigma}T_h^*$, where $\sigma$ is
as in Equation~(\ref{defofsigma}), $g(z)=(-\overline{b}z+\overline{d})^{-2}$, and
$h(z)=(cz+d)^2$ \cite{Hu}.  We follow the outline of the proof of 
Theorem~\ref{adjointcp}
to see that $\cp^*=sC_{\sigma}+K$ for some compact $K$ on $A^2$, where now
$s=[(\overline{c}\overline{\zeta}+\overline{d})/(-\overline{b}\eta+\overline{d})]^2$.  
Now the compactness of $[C_{\sigma},T_z^*]$ follows from Theorem 3 in \cite{mw},
and the compactness of $T_{\overline{h}g-\overline{h(\zeta)}g(\zeta)}C_{\sigma}$
is obtained as an application of Lemma 1 in \cite{JM} on compact Carleson
measures of the form $W(z)d(A\sigma^{-1})$, with the choice
$W(z)=|\overline{h(z)}g(z)-\overline{h(\zeta)}g(\zeta)|^2$.  We leave 
the details to
the interested reader.

\vspace{.2in}

The scalar 
$s=(\overline{c}\overline{\zeta}+\overline{d})/(-\overline{b}\eta+\overline{d})$
can equivalently be described as $|\sigma'(\eta)|$ or $|\varphi'(\zeta)|^{-1}$.  
This will
be verified below, in Proposition~\ref{sredefined}.  In particular, the 
scalar $s$ in the statement 
of Theorem~\ref{adjointcp} is strictly
positive.

\begin{cor} \label{essnorm}
For $\varphi$ a linear-fractional self-map of the disk, not an automorphism, with
$\|\varphi\|_{\infty}=1$, the self-commutator $[C_{\varphi}^*,\cp]$ is 
compact if and
only if $\varphi\circ\sigma=\sigma\circ\varphi$.
\end{cor} 
\begin{proof} We have 
$[C_{\varphi}^*,\cp]=s(C_{\varphi\circ\sigma}-C_{\sigma\circ\varphi})+K$
where $s$ is as in the statement of Theorem~\ref{adjointcp} and $K$ is compact.  
Since a difference
of non-compact linear-fractional composition operators is compact only 
if it is zero (\cite{Bo},
\cite{KM}),
the
result follows.
\end{proof}

A linear-fractional self-map whose fixed point set, relative to the Riemann sphere,
consists of a single point $\zeta$ in $\partial \D$ is termed {\it parabolic}.  
It is conjugate,
via the map $(\zeta+z)/(\zeta-z)$, to a translation by some complex number
$t$, $\mbox{Re }t\geq 0$, in the right half-plane.  When $\mbox{Re }t=0$ we have a 
(parabolic) automorphism; otherwise the map is not an automorphism.  When 
the translation
number $t$ is strictly positive, we call the associated linear-fractional self-map
of $\D$ a {\it positive parabolic} non-automorphism.  
Among the linear-fractional non-automorphisms fixing $\zeta\in\partial \D$,
the parabolic ones are characterized by $\varphi'(\zeta)=1$.  For further 
details on the classification of
linear-fractional self-maps of $\D$, see \cite{BLSS} or Chapter 0 of \cite{Sbook}.

A linear-fractional non-automorphism $\varphi$ with a fixed point 
$\zeta$ on $\partial \D$, which commutes with its Krein adjoint, 
must be parabolic.  This follows by a consideration of fixed points: 
if $\varphi$ has another
fixed point $z_0$ in the Riemann sphere, and it commutes with $\sigma$, 
then $\sigma(z_0)$
would also be fixed by $\varphi$.  Neither $\sigma(z_0)=\zeta$ nor 
$\sigma(z_0)=z_0$ are
possible, since $\sigma$ fixes the boundary point $\zeta$ if $\varphi$ 
does, and $\varphi$ fixes
$1/\overline{z_0}$ if $\sigma$ fixes $z_0$.  Thus Corollary~\ref{essnorm} 
gives another view 
of the main
result in \cite{BLSS}:  a non-automorphism linear-fractional composition 
operator $\cp$ is
non-trivially essentially normal if and only if $\varphi$ is parabolic.

\begin{prop}\label{posparabolic}
Suppose $\varphi$, not an automorphism, is a linear-fractional self-map of $\D$ 
with $\varphi(\zeta)=\eta$ for some $\zeta,\eta\in\partial \D$.  If $\sigma$ 
is the Krein
adjoint of $\varphi$, 
then $\varphi'(\zeta)\sigma'(\eta)=1$ and $\tau\equiv\varphi\circ\sigma$ 
is a positive 
parabolic non-automorphism.
\end{prop}
\begin{proof}
Using $\ \tilde{}\ $ for the Krein adjoint, we have
$\widetilde{\varphi\circ\sigma} = \tilde{\sigma}\circ\tilde{\varphi}=
\varphi\circ\sigma$.
Thus the map $\tau=\varphi\circ\sigma$, a non-automorphism fixing 
$\eta \in\partial\D$, is
its own Krein adjoint.  By the remark preceeding the statement of 
Proposition~\ref{posparabolic}, this
means that $\tau$ is parabolic and $\tau(z)=\Phi^{-1}(\Phi(z)+t)$ for 
$\Phi(z)=(\eta+z)/(\eta-z)$
and some $t$ with $\mbox{Re }t>0$.  Direct calculation, using 
$\tilde{\tau}=\tau$,
shows that $t$ must be positive. 

Since parabolic non-automorphisms have derivative one at their 
(boundary) fixed point
(\cite{Sbook}, p. 3),
we have $\varphi'(\sigma(\eta))\sigma'(\eta)=1$ or $\varphi'(\zeta)\sigma'(\eta)=1$,
as desired.
\end{proof}

The spectrum of a composition operator whose symbol is a parabolic non-automorphism 
has been described in \cite{Co83}. In particular, we have the following
result.

\begin{prop}\cite{Co83}\label{spectrum}
Let $\tau=\varphi\circ\sigma$, where $\varphi$ is a non-automorphism with
$\varphi(\zeta)=\eta$ for $\zeta,\eta\in\partial \D$.  The spectrum, $\sigma(C_{\tau})$,
and essential spectrum $\sigma_e(C_{\tau})$, are both equal to $[0,1]$.
\end{prop}

\begin{proof}
The map $\tau$ fixes $\eta\in\partial\D$, and by conjugating by a rotation, 
$C_{\tau}$ is
unitarily equivalent to a composition operator with positive parabolic 
symbol fixing $1$.
Such a map can be written as
$$\frac{(2-t)z+t}{-tz+2+t}$$ for some positive $t$.  Applying Corollary 6.2 
in \cite{Co83},
we have $\sigma(C_{\tau})=[0,1]$.  Since every point of $\sigma(C_{\tau})$ 
is a boundary
point of the spectrum, and none is isolated, we also have $\sigma_e(C_{\tau})=
\sigma(C_{\tau})=[0,1]$ (\cite{Con}, Theorem 37.8).
\end{proof}

As promised, we can describe the scalar $s$ appearing in 
Theorem~\ref{adjointcp} in a more useful
way:

\begin{prop}\label{sredefined}
Let $\varphi,\sigma$ and $s$ be as in the statement of 
Theorem~\ref{adjointcp}.  We have
$s=|\sigma'(\eta)|=|\varphi'(\zeta)|^{-1}$.
\end{prop}
\begin{proof}
Direct calculation shows that 
$$\frac{\sigma'(\eta)}{\overline{\varphi'(\zeta)}}=
\left(\frac{\overline{c}\overline{\zeta}
+\overline{d}}{-\overline{b}\eta+\overline{d}}\right)^2.$$
By Proposition~\ref{posparabolic}, $\varphi'(\zeta)=
(\sigma'(\eta))^{-1}$, so that 
$s^2=|\sigma'(\eta)|^2$.
By Theorem~\ref{adjointcp}, $C_{\varphi}C_{\varphi}^*
\equiv sC_{\varphi}C_{\sigma}\ (\mbox{mod }
\mathcal{K})=
sC_{\sigma\circ\varphi} 
$, and by Proposition~\ref{spectrum}, the essential spectrum of 
$C_{\sigma\circ\varphi}$
is $[0,1]$.  Since $C_{\varphi}C_{\varphi}^*$ is positive, the 
scalar $s$ must be positive,
and we have $s=|\sigma'(\eta)|$.
\end{proof}

\begin{cor}\label{essspec}
If $\varphi$ is a non-automorphism, linear-fractional map with 
$\varphi(\zeta)=\eta$ for
some $\zeta,\eta\in\partial\D$, then $\sigma_e(\cp^*\cp)=\sigma_e(\cp\cp^*)=
[0,s]$.
\end{cor}
\begin{proof} We have $\cp^*\equiv sC_{\sigma}\ (\mbox{mod }\mathcal{K})$ for 
$s=1/|\varphi'(\zeta)|$ by Theorem~\ref{adjointcp} and 
Proposition~\ref{sredefined}.  
Thus $\cp\cp^*\equiv sC_{\sigma\circ\varphi}\ (\mbox{mod }
\mathcal{K})$ and $\cp^*\cp\equiv sC_{\varphi\circ\sigma}\ (\mbox{mod }\mathcal{K})$, 
and the conclusion
follows from Proposition~\ref{spectrum} and 
Proposition~\ref{sredefined}.
\end{proof}

Note that since the non-zero points in $\sigma(\cp\cp^*)$ and $\sigma(\cp^*\cp)$
are the same, we also have $\sigma(\cp\cp^*)=\sigma(\cp^*\cp)$.  Moreover,
this common spectrum consists of $[0,s]$ plus at most finitely many eigenvalues
greater than $s$, and of finite multiplicity.

\section{THE UNITAL $C^*$-ALGEBRA GENERATED BY $C_{\varphi}$ AND $T_z$}

Throughout this section, $\varphi$ will be a fixed but arbitrary 
linear-fractional self-map
of $\D$ satisfying the following:
\begin{itemize}
\item[(i)] $\varphi$ is not an automorphism.
\item[(ii)] $\varphi(\zeta)=\eta$ for some $\zeta\neq\eta\in\partial\D$.
\end{itemize}
Conditions (i) and (ii) imply that $C_{\varphi}^2$ is compact on $H^2$, since
$\|\varphi\circ\varphi\|_{\infty}<1$.  

The algebra $C^*(T_z,\cp)$ is the closed linear span of all words in
$T_z, T_z^*, \cp,\cp^*$ and $I$, and contains all Toeplitz operators
$T_w$ with $w$ continuous.  
We set $\mathcal{A}=C^*(T_z,\cp)/\mathcal{K}$, and denote the cosets of $\cp$, $\cp^*$,
and $T_w$ by
$x$, $x^*$, and $t_w$, respectively.  Let $e$ denote the coset 
of the identity.
A main goal of this section is a description of $\mathcal{A}$.  
This description will
allow us, for example, to determine the essential norm and 
essential spectrum
of any element of $C^*(T_z,\cp)$.  For $\varphi$ as described 
above, $E(\varphi)=\{\zeta\}$,
and Corollary~\ref{thm1cor} implies that $T_{w-w(\zeta)}\cp$ 
is compact, that is,
$$T_w\cp\equiv w(\zeta)\cp \ \ (\mbox{mod }\mathcal{K}).$$

Since $E(\sigma)=\{\eta\}$, we also see from Corollary~\ref{thm1cor}, 
Theorem~\ref{adjointcp},
and Proposition~\ref{sredefined} that 
\begin{eqnarray*}
\cp T_w&=&(T_{\overline{w}}\cp^*)^*\\
&\equiv&s(T_{\overline{w}}C_{\sigma})^*\ \ (\mbox{mod }\mathcal{K})\\
&\equiv&s(\overline{w(\eta)}C_{\sigma})^*\ \ (\mbox{mod }\mathcal{K}),\\
&\equiv&w(\eta)\cp \ \ (\mbox{mod }\mathcal{K})
\end{eqnarray*}
where $s=|\varphi(\zeta)|^{-1}$.
In addition, $T_vT_w-T_{vw}$ is compact whenever $v$ and $w$ 
are in $C(\partial\D)$.
Phrasing these relations in terms of the cosets yields
\begin{itemize}
\item $t_wx=w(\zeta)x$
\item $xt_w=w(\eta)x$
\item $t_wx^*=w(\eta)x^*$
\item $x^*t_w=w(\zeta)x^*$
\item $t_vt_w=t_{vw}$ 
\end{itemize}
for all $w$ and $v$ in $C(\partial\D)$.
Since $x^2=(x^*)^2=0$, we generate $\mathcal{A}$ as a Banach space
from linear combinations of 
$$t_w,\ (x^*x)^m,\ (xx^*)^n,\ x(x^*x)^j,\ x^*(xx^*)^k,$$ where 
$w\in C(\partial\D)$, the 
integers 
$m,n$ are positive, and the integers $j$ and $k$ are non-negative. 

Let $K$ be a compact subset of the non-negative real 
numbers which contains
$[0,s]$.  We write $C_0(K)$ for the space of functions in $C(K)$ which vanish
at zero.  
We will need the next result, which follows easily from
the Hahn-Banach theorem and the Riesz Representation Theorem; here $t$
denotes the independent variable.

\begin{lemma}\label{hbrrt}
\begin{itemize}
\item[(i)] Let $\mathcal{R}$ and $\mathcal{S}$ be 
dense linear manifolds in $C_0(K)$
and $C(K)$, respectively.  If $\alpha>0$, then
$$\overline{t^{\alpha}\mathcal{R}}=\overline{t^{\alpha}\mathcal{S}}=C_0(K).$$
\item[(ii)] Suppose $0<\lambda\leq s$ and let $\mathcal{T}$ be a 
linear manifold which is dense in the subspace $\{f\in C(K):f(\lambda)=0\}$.  
Then
$$\overline{t^{\alpha}\mathcal{T}}=\{f\in C_0(K):f(\lambda)=0\}.$$
\end{itemize}
\end{lemma}

We next introduce the various objects which are central to 
our analysis and
record some observations about them.

\subsection{THE $C^*$-ALGEBRA $\mathcal{C}$}  
It follows from the relations described above that for every 
continuous function
$w$ on $\partial\D$, $t_w$ commutes with $x^*x$ and $xx^*$.  
Further, if we let $C_{\zeta,\eta}$
denote the algebra of all $w$ in $C(\partial \D)$ 
satisfying $w(\eta)=w(\zeta)$,
then $t_w$ commutes with $x$ and $x^*$ whenever $w$ 
lies in $C_{\zeta,\eta}(\partial \D)$.  
Finally note that the 
self-adjoint element $a\equiv xx^*+x^*x$ commutes with 
both $x$ and
$x^*$.
The spectrum of $a$ is easily identified:

\begin{prop}\label{speca}
Let $x$ be the coset of $\cp$ in $\mathcal{A}$, 
where $\varphi=(az+b)/(cz+d)$ 
satisfies conditions (i)-(ii) stated at the beginning of Section 4.  
If $a=xx^*+x^*x$, then $\sigma(a)=\sigma(xx^*)\cup
\sigma(x^*x)=[0,s]$ where $s=1/|\varphi'(\zeta)|$.
\end{prop}

\begin{proof}
The elements $x^*x$ and $xx^*$ generate a commutative $C^*$-algebra.
It follows from Gelfand theory, the facts that 
$(xx^*)(x^*x)=(x^*x)(xx^*)=0$, and (by
Corollary~\ref{essspec}) $\sigma(xx^*)
=\sigma(x^*x)=[0,s]$, that
$\sigma(a)=\sigma(xx^*)\cup\sigma(x^*x)$.
\end{proof}

Let $\mathcal{C}$ denote the (necessarily commutative) $C^*$-algebra generated by 
$a$ and the Toeplitz cosets $\{t_w:w\in C_{\zeta,\eta}(\partial\D)\}$.  
Clearly, $\mathcal{C}$ lies in
the center of $\mathcal{A}$.  We next describe the Gelfand theory of $\mathcal{C}$.
First we look at the algebra $C_{\zeta,\eta}(\partial \D)$.  

It is easy to see that the multiplicative linear functionals on 
$C_{\zeta,\eta}(\partial\D)$
are all point evaluations
$$\ell_{\lambda}:f\rightarrow f(\lambda)$$
with the proviso that $\ell_{\eta}=\ell_{\zeta}$.  Accordingly, the maximal ideal
space of $C_{\zeta,\eta}(\partial\D)$ is a ``figure eight", namely, 
the circle $\partial \D$
with $\zeta$ and $\eta$ identified.  
We denote by $\Lambda$ the disjoint union of $\partial\D$ and $[0,s]$, 
with $\zeta,\eta$
and $0$ identified to a point ${\bf p}$ (a figure eight with an interval attached).  
Given $w$ in $C_{\zeta,\eta}(\partial\D)$, let us agree to extend $w$ continuously
to $\Lambda$ by setting $w(\lambda)=w(\zeta)=w(\eta)$ when $\lambda ={\bf p}$
or $0<\lambda \leq s$.  Similarly, if $f\in C_0([0,s])$, extend $f$ continuously
to $\Lambda$ by putting $f({\bf p})=f(0)=0$ and $f(\lambda)=0$ for
$\lambda\in\partial\D\backslash\{\zeta,\eta\}$.  With these understandings, 
which remain
in force throughout, we have the following result.

\begin{prop}\label{descC}
The algebra $\mathcal{C}$ consists of all elements of the form
$b=t_w+f(a)$ where $w$ is in $C_{\zeta,\eta}(\partial\D)$ and $f$ is 
in $C_0([0,s])$.
Moreover, $b$ uniquely determines $w$ and $f$.  The maximal ideal space 
of $\mathcal{C}$ 
coincides with $\Lambda$, and the Gelfand transform from $\mathcal{C}$ 
to $C(\Lambda)$
has the form
$$t_w+f(a)\rightarrow w+f.$$
\end{prop}

\begin{proof} We temporarily write $\mathcal{C}_0$ for 
$\{t_w+f(a):w\in C_{\zeta,\eta}(\partial\D)
\mbox{ and }f\in C_0([0,s])\}$.  
If $w(\zeta)=w(\eta)$ and $f$ is in $C_0([0,s])$, then,
since $f$ is a uniform limit of polynomials vanishing at zero
(and $(x^*x)(xx^*)=0$), we have
\begin{eqnarray*}
t_wf(a)&=&t_w(f(x^*x)+f(xx^*))\\
&=&w(\eta)f(x^*x)+w(\zeta)f(xx^*)\\&=&w(\zeta)f(a).
\end{eqnarray*}
Since $t_wt_v=t_{wv}$ for continuous $w$ and $v$, we see that
$\mathcal{C}_0$ is an algebra.  

Suppose $\ell$ is a multiplicative linear functional on $\mathcal{C}$.  
Restricting $\ell$ to
$$\{t_w:w\in C_{\zeta,\eta}(\partial \D)\}
\cong C_{\zeta,\eta}(\partial \D)$$
we see that there is a unique $\alpha\in\partial\D$ with $\ell(t_w)=w(\alpha)$
for all continuous $w$ with $w(\zeta)=w(\eta)$.  Restricting $\ell$ to
$$\{f(a):f\in C([0,s])\}\cong C([0,s])$$
shows that there is a unique point $\beta$ in $[0,s]$ with
$\ell(f(a))=f(\beta)$.  Thus
$$\ell(t_wf(a))=\ell(t_w)\ell(f(a))=w(\alpha)f(\beta).$$
Also, if $f(0)=0$, then $t_wf(a)=w(\zeta)f(a)$ as seen above,
so $\ell(t_wf(a))=w(\zeta)f(\beta)$.  Since any function in
$C_0([0,s])$ vanishes at $0$, we can have $\alpha\in \partial\D\backslash
\{\zeta,\eta\}$ if $\beta=0$, but if $0<\beta\leq s$, $\alpha\in\{\zeta,\eta\}$.
Thus with the understandings stated prior to the statement of the
proposition, $\ell(t_w+f(a))=w(\lambda)+f(\lambda)$ for
a unique $\lambda$ in $\Lambda$ and any $t_w+f(a)$ in $\mathcal{C}_0$.

The above arguments show that $C(\Lambda)$ is the Gelfand representation
for $\mathcal{C}$.  Moreover, the map
$$t_w+f(a)\rightarrow w+f$$
from $\mathcal{C}_0$ to $C(\Lambda)$ is an isometric $*$- homomorphism
from $\mathcal{C}_0$ to  $C(\Lambda)$.  But $C(\Lambda)$ consists of
exactly such sums $w+f$, so this $*$-homomorphism is onto $C(\Lambda)$.  
Since $C(\Lambda)$ is complete, so is $\mathcal{C}_0$.
Since $\mathcal{C}_0$ is dense in $\mathcal{C}$, we conclude
$\mathcal{C}_0=\mathcal{C}$.
\end{proof}

\subsection{THE POLAR DECOMPOSITION OF $C_{\varphi}$ AND THE ALGEBRA 
$\mathcal{A}_0$}
We 
begin with some observations on the polar decomposition of any operator $T$
on a Hilbert space $\mathcal{H}$.  
Suppose that $T=U\sqrt{T^*T}$, where $U$ is a partial isometry with
initial space $(\mbox{ker }T)^{\perp}=\overline{T^*\mathcal{H}}$
and final space $\overline{T\mathcal{H}}=(\mbox{ker }T^*)^{\perp}$.
The operators $U^*U$ and $UU^*$ are, respectively, the projections
onto $(\mbox{ker }T)^{\perp}$ and $\overline{T\mathcal{H}}$.  
Moreover, $UT^*T=TT^*U$ and so 
\begin{equation}\label{halmos}
Uf(T^*T)=f(TT^*)U
\end{equation} for 
all functions continuous on the spectra of both $T^*T$ and $TT^*$.  Taking
$f$ to be the square root function shows that the polar decomposition
for $T^*$ is $T^*=U^*\sqrt{TT^*}$.
The partial
isometry $U$ is unitary if $T$ and $T^*$ are one-to-one.  Observe that
every non-trivial composition operator is one-to-one, and the adjoint formula of
Equation~(\ref{cowenadjoint}) guarantees that, for linear-fractional 
composition operators,
the adjoint is also one-to-one.  Thus the linear-fractional 
composition operators
under consideration here have the polar decomposition 
$\cp=U\sqrt{\cp^*\cp}$ where
$U$ is {\it unitary}.
If we apply these remarks to $T=C_{\varphi}=U\sqrt{C_{\varphi}^*C_{\varphi}}$,
we have
$x=u\sqrt{x^*x}$ and $x^*=u^*\sqrt{xx^*}$
where $u=[U]$, the coset of $U$ modulo $\mathcal{K}$, and $x=[C_{\varphi}]$.  
Moreover, as observed above,
$U$, and hence $u$, are unitary.
By Corollary~\ref{essspec}, the sets $\sigma(x^*x)=\sigma_e(\cp^*\cp)$ and
$\sigma(xx^*)=\sigma_e(\cp\cp^*)$ both coincide with $[0,s]$,
where $s=|\varphi'(\zeta)|^{-1}$.

Now $C^*(T_z,\cp)$ is the closed linear span of elements of the form
$$T_w,\  f(\cp^*\cp), \ g(\cp\cp^*),\ \cp p(\cp^*\cp), \ \cp^*q(\cp\cp^*),\ K,$$  
where $f,g,p$ and $q$ are polynomials
with $f(0)=g(0)=0$, $w$ is in $C(\partial\D)$, and $K$ is a 
compact operator.  
The map $f\rightarrow f(\cp^*\cp)$ extends to a $*$-
isomorphism of $C_0(\sigma(\cp^*\cp))$ onto the closed subspace
$\{f(\cp^*\cp):f\in C_0(\sigma(\cp^*\cp))\}$ in $\mathcal{B}(H^2)$; the
analogous statement holds for the map $g\rightarrow g(\cp\cp^*)$.  
Writing 
$$\cp p(\cp^*\cp)=U\sqrt{\cp^*\cp}p(\cp^*\cp),$$
we see by Lemma~\ref{hbrrt}
that 
$$\overline{\{\cp p(\cp^*\cp): p\mbox{ a polynomial}\}}= 
\{Uh(\cp^*\cp):h\in C_0(\sigma(\cp^*\cp))\};$$ similarly, 
$$\overline{\{\cp^* q(\cp\cp^*): q\mbox{ a polynomial}\}}=
\{U^*k(\cp\cp^*):k\in C_0(\sigma(\cp\cp^*))\}.
$$  Thus $\mathcal{A}=C^*(T_z,\cp)/\mathcal{K}$ contains,
and is the closure of, the set $\mathcal{A}_0$ of elements of the form
\begin{equation}\label{defofa0}
b=t_w+f(x^*x)+g(xx^*)+uh(x^*x)+u^*k(xx^*)
\end{equation}
where $w\in C(\partial \D)$, and $f,g,h$ and $k$ 
are in $C_0([0,s])$, with $s=1/|\varphi'(\zeta)|$.
We will see later that $\mathcal{A}_0=\mathcal{A}$; for now we
show that $\mathcal{A}_0$ is an algebra, and
each element of $\mathcal{A}_0$ has a unique representation in
the above form.  To this end,
we record some consequences of the next pair of equations, which follow
from Equation~(\ref{halmos}) by taking cosets and adjoints:
\begin{equation}\label{cosethalmos}
uf(x^*x)=f(xx^*)u \mbox{ and } u^*f(xx^*)=f(x^*x)u^*
\end{equation}
for all $f\in C([0,s])$.

\begin{prop}\label{a0analg}
If $\mathcal{A}_0$ is defined as above, then $\mathcal{A}_0$ is 
an algebra.
\end{prop}
\begin{proof}
We must show that given elements $b_1\in\mathcal{A}_0$ and 
$b_2\in\mathcal{A}_0$ having the
form 
$$b_j=t_{w_j}+f_j(x^*x)+g_j(xx^*)+uh_j(x^*x)+u^*k_j(xx^*),\  j=1,2$$
with $w_j\in C(\partial\D)$ and $f_j,g_j,h_j,k_j$ in $C_0([0,s])$,
then $b_1b_2$ has the
same form.  To do this, it suffices to show that that the product
of any of the five terms of $b_1$ with any of the five terms of $b_2$
is again in $\mathcal{A}_0$.  Some of these verifications are immediate,
for example $f_1(x^*x)f_2(x^*x)=f_1f_2(x^*x)$, where $f_1f_2$
is in $C_0([0,s])$ if $f_1$ and $f_2$ are.  
For the others, we make use of the basic equations of (\ref{cosethalmos}) together
with:
\begin{equation}\label{productzero}
f(x^*x)g(xx^*)=0=g(xx^*)f(x^*x)
\end{equation}
for $f$ and $g$
in $C_0([0,s])$.
Equation~(\ref{productzero}) follows by uniformly approximating 
$f$ and $g$ by polynomials vanishing
at $0$.  From these equations we see that
\begin{itemize}
\item $g_1(xx^*)uh_2(x^*x)=ug_1(x^*x)h_2(x^*x)$,
\item $uh_1(x^*x)g_2(xx^*)=0$, 
\item $uh_1(x^*x)uh_2(x^*x)=uh_1(x^*x)h_2(xx^*)u^*=0$,
\item $uh_1(x^*x)u^*k_2(xx^*)=h_1(xx^*)uu^*k_2(xx^*)=h_1(xx^*)k_2(xx^*)$,
\item $u^*k_1(xx^*)uh_2(x^*x)=u^*uk_1(x^*x)h_2(x^*x)=k_1(x^*x)h_2(x^*x)$,
and
\item$u^*k_1(xx^*)u^*k_2(xx^*)=u^*k_1(xx^*)k_2(x^*x)u^*=0.$
\end{itemize}
Similarly we see (using the coset identities preceeding Lemma~\ref{hbrrt})
that for $f,g,h$, and $k$ in $C_0([0,s])$ and $w\in C(\partial \D)$, 
\begin{itemize}
\item $t_wf(x^*x)=w(\eta)f(x^*x)$,
\item $t_wg(xx^*)=w(\zeta)g(xx^*)$,
\item $t_wuh(x^*x)=w(\zeta)uh(x^*x)$,
\item $t_wu^*k(xx^*)=w(\eta)u^*k(xx^*)$.
\end{itemize}
This shows that $\mathcal{A}_0$ is an algebra.
\end{proof}

The next result addresses the uniqueness of representation of 
elements in $\mathcal{A}_0$.
\begin{prop}
For an element $b$ in $\mathcal{A}_0$, there is a unique 
$w\in C(\partial \D)$  and
unique functions $f,g,h$ and $k$ in $C_0([0,s])$ so that 
Equation~(\ref{defofa0}) holds.
\end{prop}
\begin{proof}
It suffices to show that if 
\begin{equation}\label{uniqueness}
0=t_w+f(x^*x)+g(xx^*)+uh(x^*x)+u^*k(xx^*),
\end{equation}
then each term on the right-hand side is zero.  Multiplying on the 
right by $x^*x$ yields
\begin{eqnarray*}
0&=&t_wx^*x+f(x^*x)x^*x+g(xx^*)x^*x+uh(x^*x)x^*x+u^*k(xx^*)x^*x\\
&=&w(\eta)x^*x+f(x^*x)x^*x+uh(x^*x)x^*x
\end{eqnarray*}
so that
$$
uh(x^*x)x^*x=-[w(\eta)x^*x+f(x^*x)x^*x].$$
The right-hand side is normal, and the left-hand side has square zero, 
so both sides
must vanish.  Thus $h\equiv0$ and $f+w(\eta)\equiv 0$ on $[0,s]$; 
since $f(0)=0$, we
must have $w(\eta)=0$
and $f\equiv 0$.  
Thus Equation~(\ref{uniqueness}) is now
$$0=t_w+g(xx^*)+u^*k(xx^*).$$
Multiplying on the left by $xx^*$ gives
\begin{eqnarray*}
0&=&xx^*t_w+xx^*g(xx^*)+xx^*u^*k(xx^*)\\
&=&w(\zeta)xx^*+xx^*g(xx^*)+xx^*u^*k(xx^*)
\end{eqnarray*}
so that
$$-[w(\zeta)xx^*+xx^*g(xx^*)]=xx^*u^*k(xx^*)=0.$$
It follows that $g+w(\zeta)\equiv 0$ on $[0,s]$; since
$g(0)=0$, we see that $w(\zeta)=0$ and $g\equiv 0$ on $[0,s]$.  Returning again
to Equation~(\ref{uniqueness}) we have
$$0=t_w+u^*k(xx^*).$$
Multiplying on the left by $x^*x$ yields
\begin{eqnarray*}
0&=&x^*xt_w+x^*xu^*k(xx^*)\\
&=&w(\eta)x^*x+x^*xk(x^*x)u^*.
\end{eqnarray*}
Since $w(\eta)=0$, this forces $k\equiv 0$, and from this
it follows finally that $t_w=0$.

\end{proof}

\subsection{LOCALIZATION AND THE STRUCTURE OF $\mathcal{A}$}
For $\lambda$ in $\Lambda$, let $\mathcal{I}_{\lambda}$ denote the closed,
two-sided ideal in $\mathcal{A}$ generated by the maximal ideal
$$J_{\lambda}=\{t_w+f(a):w\in C_{\zeta,\eta}(\partial\D),f
\in C_0([0,s]) \mbox{ and }w(\lambda)+f(\lambda)=0\}$$
of $\mathcal{C}$.  Here $w$ and $f$ are understood to extend
to $\Lambda$ as described prior to Proposition~\ref{descC}.  For 
$b$ in $\mathcal{A}$, we 
write $[b]_{\mathcal{I}_{\lambda}}$ for the coset of $b$ in 
$\mathcal{A}/\mathcal{I}_{\lambda}$.
The localization theorem of R. G. Douglas (\cite{Do}, p. 196) tells 
us that
$$\|b\|=\sup_{\lambda\in \Lambda}\|[b]_{\mathcal{I}_{\lambda}}\|,$$
and the map
$$b\rightarrow \{[b]_{\mathcal{I}_{\lambda}}\}_{\lambda\in\Lambda}$$ is
an isometric $*$- homomorphism of $\mathcal{A}$ into
$\sum_{\lambda\in\Lambda}\oplus \mathcal{A}/\mathcal{I}_{\lambda}$.  Moreover,
a given $b$ in $\mathcal{A}$ is invertible if and only if each
coset $[b]_{\lambda}$ is invertible, for $\lambda \in \Lambda$.  Our immediate 
objective is to
compute the local algebras $\mathcal{A}/\mathcal{I}_{\lambda}$.

For $\lambda$ in $\Lambda$ we define a map $\Phi_{\lambda}:\mathcal{A}_0
\rightarrow \mathbb{M}_2$, the algebra of $2\times 2$ matrices, as follows.  Let
$b$ in $\mathcal{A}_0$ be given by Equation~(\ref{defofa0}).  We put 
\begin{equation}\label{defphil}
\Phi_{\lambda}(b)=
\left\{\begin{array}{l}
\left[\begin{array}{lr}
w(\zeta)+g(\lambda)&h(\lambda)\\
k(\lambda)&w(\eta)+f(\lambda)
\end{array}\right] \ \mbox{ if }0<\lambda\leq s,
\\[1mm]
\\
\left[\begin{array}{lr}
w(\zeta)&0\\
0&w(\eta)
\end{array}\right]\ \mbox {if }\lambda={\bf p},
\\[1mm]
\\
\left[\begin{array}{lr}w(\lambda)&0\\
0&w(\lambda)
\end{array}\right]\ \mbox{ if }\lambda \in \partial \D\backslash\{\zeta,\eta\}.
\end{array}
\right.
\end{equation}
We write $I_{2\times 2}$ for the identity matrix in 
$\mathbb{M}_2$ and $\mathbb{M}_2^{diag}$ for the algebra
of $2\times 2$ diagonal matrices.  The range of $\Phi_{\lambda}$ will be denoted
$\mbox{Ran }\Phi_{\lambda}$.  

\begin{prop}\label{rangephi}
For each $\lambda$ in $\Lambda$, $\Phi_{\lambda}$ is a $*$- homomorphism
from $\mathcal{A}_0$ to $\mathbb{M}_2$ with
\begin{equation}
\mbox{Ran }\Phi_{\lambda}=
\left\{\begin{array}{l}
\mathbb{M}_2 \mbox{ when }0<\lambda\leq s,\\
\mathbb{M}_2^{diag}\mbox{ when }\lambda = {\bf p},\\
\{cI_{2\times 2}:c\in\mathbb{C}\} 
\mbox{ when }\lambda\in\partial\D\backslash\{\zeta,\eta\}.
\end{array}
\right.
\end{equation}
\end{prop}

\begin{proof}
First consider $\lambda >0$.  Any element $b$ in $\mathcal{A}_0$ has 
the form
$b=t_w+y$, where $w$ is in $C(\partial\D)$ and
\begin{equation}\label{ylocal}
y=f(x^*x)+g(xx^*)+uh(x^*x)+u^*k(xx^*)
\end{equation}
with $f,g,h,k$ in $C_0([0,s])$.  Given $b_1=t_{w_1}+y_1$ and
$b_2=t_{w_2}+y_2$ in $\mathcal{A}_0$, 
\begin{equation}\label{b1b2}
b_1b_2=t_{w_1}t_{w_2}+y_1t_{w_2}+t_{w_1}y_2+y_1y_2
\end{equation}
Taking the notation from Equation~(\ref{ylocal}) for $y_1$ and
$y_2$, we have 
\begin{eqnarray*}
y_1y_2&=&[f_1(x^*x)f_2(x^*x)+k_1(x^*x)h_2(x^*x)]\\&+&
u[g_1(x^*x)h_2(x^*x)+h_1(x^*x)f_2(x^*x)]\\
&+&u^*[k_1(xx^*)g_2(xx^*)+f_1(xx^*)k_2(xx^*)]\\&+&
[g_1(xx^*)g_2(xx^*)+h_1(xx^*)k_2(xx^*)],
\end{eqnarray*}
where we have used the list of identities in the proof of 
Proposition~\ref{a0analg} and collected like terms.
Thus
\begin{eqnarray*}
\Phi_{\lambda}(y_1y_2)&=&
\left[\begin{array}{lr}
g_1(\lambda)g_2(\lambda)+h_1(\lambda)k_2(\lambda)&
g_1(\lambda)h_2(\lambda)+h_1(\lambda)f_2(\lambda)\\
k_1(\lambda)g_2(\lambda)+f_1(\lambda)k_2(\lambda)&
f_1(\lambda)f_2(\lambda)+k_1(\lambda)h_2(\lambda)
\end{array}\right]\\
&=&\Phi_{\lambda}(y_1)\Phi_{\lambda}(y_2).
\end{eqnarray*}

Now 
$$t_{w_1}y_2=w_1(\eta)f_2(x^*x)+w_1(\zeta)g_2(xx^*)+
w_1(\zeta)uh_2(x^*x)+w_1(\eta)u^*k_2(xx^*).$$
Thus 
\begin{eqnarray*}
\Phi_{\lambda}(t_{w_1}y_2)&=&
\left[\begin{array}{lr}
w_1(\zeta)g_2(\lambda)&w_1(\zeta)h_2(\lambda)\\
w_1(\eta)k_2(\lambda)&w_1(\eta)f_2(\lambda)
\end{array}\right]\\
&=&\left[\begin{array}{lr}
w_1(\zeta)&0\\
0&w_1(\eta)
\end{array}\right]
\left[\begin{array}{lr}
g_2(\lambda)&h_2(\lambda)\\
k_2(\lambda)&f_2(\lambda)
\end{array}\right]\\
&=&\Phi_{\lambda}(t_{w_1})\Phi_{\lambda}(y_2).
\end{eqnarray*}

Similarly, we find $\Phi_{\lambda}(y_1t_{w_2})
=\Phi_{\lambda}(y_1)\Phi_{\lambda}(t_{w_2})$.
Since $t_{w_1}t_{w_2}=t_{w_1w_2}$, it follows that
$$
\Phi_{\lambda}(t_{w_1}t_{w_2})=
\Phi_{\lambda}(t_{w_1})\Phi_{\lambda}(t_{w_2}).
$$
Applying $\Phi_{\lambda}$ to both sides of Equation~(\ref{b1b2}) and 
invoking the above identities, we see that
\begin{eqnarray*}
\Phi_{\lambda}(b_1b_2) &=&\Phi_{\lambda}(t_{w_1})\Phi_{\lambda}(t_{w_2})
+\Phi_{\lambda}(y_1)\Phi_{\lambda}(t_{w_2})\\
&+&\Phi_{\lambda}(t_{w_1})
\Phi_{\lambda}(y_2)+\Phi_{\lambda}(y_1)\Phi_{\lambda}(y_2)\\
&=&\left(\Phi_{\lambda}(t_{w_1})+\Phi_{\lambda}(y_1)\right)\left(
\Phi_{\lambda}(t_{w_2})+\Phi_{\lambda}(y_2)\right)\\
&=&\Phi_{\lambda}(b_1)\Phi_{\lambda}(b_2)
\end{eqnarray*}
as desired.  Clearly the range of $\Phi_{\lambda}$ is $\mathbb{M}_2$,
which yields the conclusion for $0<\lambda\leq s$.  

The remaining cases $\lambda={\bf p}$ and $\lambda\in \partial\D\backslash\{\zeta,
\eta\}$, which are considerably easier since there one has $\Phi_{\lambda}(t_w+y)=
\Phi_{\lambda}(t_w)$, are left for the reader.
\end{proof}

\begin{prop}\label{kerphi}
For $\lambda \in \Lambda$, $\overline{\mbox{ker }\Phi_{\lambda}}
=\mathcal{I}_{\lambda}$.  
\end{prop}
\begin{proof} 
For $\lambda$ in $\Lambda$, denote by $\mathcal{I}_{\lambda}^{alg}$
the two-sided algebraic ideal in $\mathcal{A}_0$ generated by $J_{\lambda}$.  
Since $\mbox{ker }\Phi_{\lambda}$ is an ideal containing
$J_{\lambda}$, we know
$$J_{\lambda}\subset \mathcal{I}_{\lambda}^{alg}\subset
\mbox{ker }\Phi_{\lambda}.$$
By definition, $\mathcal{I}_{\lambda}=\overline{\mathcal{I}_{\lambda}^{alg}}$.
It suffices to show that 
$\mbox{ker }\Phi_{\lambda}\subset \overline{\mathcal{I}_{\lambda}
^{alg}}$, for then we will have
$$\mathcal{I}_{\lambda}=
\overline{\mathcal{I}_{\lambda}^{alg}}\subset\overline{\mbox{ker }\Phi_{\lambda}}
\subset \overline{\mathcal{I}_{\lambda}^{alg}}=\mathcal{I}_{\lambda},$$
which gives the desired conclusion.  

Consider first the case $0<\lambda\leq s$.  An element $b$ in $\mathcal{A}_0$, given
by Equation~(\ref{defofa0}), lies in $\mbox{ker }\Phi_{\lambda}$ exactly when 
$w(\zeta)+g(\lambda), w(\eta)+f(\lambda), h(\lambda)$ and $k(\lambda)$ are all
zero.  We claim that the sum of the first three terms on the right side of 
Equation~(\ref{defofa0})
lie in $\mathcal{I}_{\lambda}^{alg}$.  To see this, pick $m$ and $n$ in $C(\partial\D)$
with $m+n\equiv 1, m(\zeta)=0, m(\eta)=1$, and $n(\zeta)=1, n(\eta)=0$.  Then
$w=mw+nw$ so that $t_w=t_{mw}+t_{nw}$.  To prove the
claim, it is enough to show that both $t_{mw}+f(x^*x)$ and $t_{nw}+g(xx^*)$
lie in $\mathcal{I}_{\lambda}^{alg}$.  Consider $t_{mw}+f(x^*x)$  

{\bf Case 1: $w(\eta)\neq 0$.}  Putting $m_1=mw/w(\eta)$, we see
that
\begin{eqnarray*}
t_{mw}+f(x^*x)&=&t_{mw}+m_1(\eta)f(x^*x)+m_1(\zeta)f(xx^*)\\
&=&t_{m_1w(\eta)}+t_{m_1}(f(x^*x)+f(xx^*))\\
&=&t_{m_1}(t_{w(\eta)}+f(a)).
\end{eqnarray*}
Since $w(\eta)$ is constant (and hence lying in $C_{\zeta,\eta}(\partial\D)$)
and $w(\eta)+f(\lambda)=0$, $t_{w(\eta)}+f(a)$ lies in $J_{\lambda}$, so 
$t_{mw}+f(x^*x)\in \mathcal{I}_{\lambda}^{alg}$.

{\bf Case 2: $w(\eta)=0$}.  If $m$ and $n$ are as above, $mw$ vanishes at
both $\zeta$ and $\eta$.  Fix a closed arc $I$ in $\partial\D$ whose
interior contains $\zeta$, but with $\eta$ not in $I$.  This time,
define $m_1=|mw|^{1/2}$ on $I$, $m_1>0$ on $\partial\D\backslash I$,
and $m_1(\eta)=1$.  Let $w_1$ be $mw/|mw|^{1/2}$ when $mw\neq 0$ and
$0$ otherwise.  Note that $w_1$ is continuous and $mw=m_1w_1$ on $\partial\D$.
Thus
\begin{eqnarray*}
t_{mw}+f(x^*x)&=& t_{mw}+m_1(\eta)f(x^*x)+m_1(\zeta)f(xx^*)\\
&=&t_{m_1w_1}+t_{m_1}(f(x^*x)+f(xx^*))\\
&=&t_{m_1}(t_{w_1}+f(a)).
\end{eqnarray*}
Since $w_1(\zeta)=w_1(\eta)=0$ and $f(\lambda)=0$, 
$t_{w_1}+f(a)$ lies in $J_{\lambda}$.  We conclude
that $t_{mw}+f(x^*x)$ is in $\mathcal{I}_{\lambda}^{alg}$ in Case 2,
as well as Case 1.  A similar argument shows that $t_{nw}+g(xx^*)$
lies in $\mathcal{I}_{\lambda}^{alg}$ in both cases, thus proving the claim.  

Next we show that the fourth term in $b$, $uh(x^*x)$, is in $\overline{
\mathcal{I}_{\lambda}^{alg}}$.  If $p$ is continuous on $[0,s]$, with
$p(0)=p(\lambda)=0$, then $p(a)$ lies in $J_{\lambda}$.  Thus $
xp(x^*x)=xp(a)$ is in $\mathcal{I}_{\lambda}^{alg}$.  Writing
$x=u\sqrt{x^*x}$, we see that $xp(a)=u\sqrt{x^*x}p(x^*x)$.
According to (ii) of Lemma~\ref{hbrrt}, the closure of 
such objects includes
our fourth term $uh(x^*x)$, so that $uh(x^*x)$ is in $\overline{\mathcal{I}_{\lambda}
^{alg}}$.  
Similarly, $\overline{\mathcal{I}_{\lambda}^{alg}}$ 
contains $u^*k(xx^*)$,
the fifth term of $b$, so that $b$ is in 
$\overline{\mathcal{I}_{\lambda}^{alg}}$
as desired.  This completes the proof for $0<\lambda\leq s$.  

Next we consider the case $\lambda ={\bf p}= \{0,\zeta,\eta\}$, the
triple point in $\Lambda$.  Recall that if $f$ is in $C_0([0,s])$, then
$f({\bf p})=f(0)=0$, while any $w$ in $C_{\zeta,\eta}(\partial\D)$ satisfies
$w({\bf p})=w(\zeta)=w(\eta)$.  An element $b$ of $\mathcal{A}_0$, specified
by Equation~(\ref{defofa0}), lies in the kernel of $\Phi_{{\bf p}}$ exactly
when $w(\zeta)=w(\eta)=0$.  We want to show that $\mbox{ker }\Phi_{{\bf p}}
\subset \overline{\mathcal{I}_{{\bf p}}^{alg}}$.  Let $m$ and $n$ be as 
described above.  For $f$ in $C_0([0,s])$, 
\begin{eqnarray*}
t_mf(a)&=&t_m(f(x^*x)+f(xx^*))\\
&=&m(\eta)f(x^*x)+m(\zeta)f(xx^*)\\
&=&f(x^*x),
\end{eqnarray*}
and similarly, for $g\in C_0([0,s])$, 
$t_ng(a)=g(xx^*)$.  
Thus $f(x^*x)$ and $g(xx^*)$ lie in $\mathcal{I}_{{\bf p}}^{alg}$.  If
$w(\zeta)=w(\eta)=0$, then $t_w$ lies in $J_{{\bf p}}\subset
 \mathcal{I}_{{\bf p}}^{alg}$.
As noted above for the case $0<\lambda\leq s$, $uh(x^*x)$ and 
$u^*k(xx^*)$ both lie in
$\overline{\mathcal{I}_{\lambda}^{alg}}$ and thus so does $b$, establishing the 
conclusion for $\lambda = {\bf p}$.  

Finally, if $\lambda$ is in $\partial \D\backslash\{\zeta,\eta\}$, note
that $J_{\lambda}$ consists of those elements $t_w+f(a)$ with $w(\lambda)=0$,
while the elements of $\mbox{ker }\Phi_{\lambda}$ have the form given by
Equation~(\ref{defofa0}), with $w(\lambda)=0$.  It follows easily 
(and similarly), that
$\overline{\mathcal{I}_{\lambda}^{alg}}$ contains $\mbox{ker }\Phi_{\lambda}$ in this
case as well.
\end{proof}

\begin{prop}\label{quotient}
Let $\lambda \in \Lambda$.  
\begin{itemize}
\item [(i)] If $0<\lambda\leq s, \mathcal{A}/\mathcal{I}_{\lambda}$ is 
$*$-isomorphic
to $\mathbb{M}_2$.  
\item[(ii)] $\mathcal{A}/I_{{\bf p}}$ is $*$-isomorphic to $\mathbb{M}_2^{diag}$.
\item[(iii)] If $\lambda$ is in $\partial\D\backslash\{\zeta,\eta\}$, $\mathcal{A}/
\mathcal{I}_{\lambda}$ is $*$-isomorphic to $\{cI_{2\times 2}:c\in\mathbb{C}\}$.
\end{itemize}
\end{prop}
\begin{proof} 
For an ideal $\mathcal{I}$ in an algebra $\mathcal{B}$, we 
write $[b]_{\mathcal{I}}$ throughout
for the coset in $\mathcal{B}/\mathcal{I}$ of an element $b$ in $\mathcal{B}$.  
First suppose $0<\lambda\leq s$.  Since $\mbox{ker }
\Phi_{\lambda}\subset \mathcal{A}_0\cap \mathcal{I}_{\lambda}$,
we may define a $*$-homomorphism $$\Gamma_{\lambda}:
\mathcal{A}_0/\mbox{ker }\Phi_{\lambda}
\rightarrow \mathcal{A}_0/(\mathcal{A}_0\cap \mathcal{I}_{\lambda})$$ by
$$\Gamma_{\lambda}([b]_{\mbox{ker }\Phi_{\lambda}})=
[b]_{(\mathcal{A}_0\cap \mathcal{I}_{\lambda})}.$$
By Proposition~\ref{rangephi} we know that $\mathcal{A}_0/\mbox{ker }\Phi_{\lambda}$ 
is $^*$-isomorphic
to $\mathbb{M}_2$; write this isomorphism as $T_{\lambda}:\mathbb{M}_2
\rightarrow \mathcal{A}_0/\mbox{ker }\Phi_{\lambda}$.  Thus we have a sequence 
of onto $*$-
homomorphisms

\begin{equation}\label{3hom}
\mathbb{M}_2\rightarrow\mathcal{A}_0/\mbox{ker }\Phi_{\lambda}\rightarrow
\mathcal{A}_0/(\mathcal{A}_0\cap \mathcal{I}_{\lambda})
\rightarrow(\mathcal{A}_0+\mathcal{I}_{\lambda})/\mathcal{I}_{\lambda},
\end{equation}
where the first map is $T_{\lambda}$, the second is
$\Gamma_{\lambda}$ and the last, call it $R_{\lambda}$, is provided by the 
first isomorphism theorem for rings 
(see, for example,
p. 105 in \cite{Ja})
and has the form 
$$R_{\lambda}:[b]_{\mathcal{A}_0\cap \mathcal{I}_{\lambda}}
\rightarrow [b]_{\mathcal{I}_{\lambda}}.$$
Since $\mathcal{A}_0$ is dense in $\mathcal{A}$, 
so is
$\mathcal{A}_0+\mathcal{I}_{\lambda}$, and we have 
$(\mathcal{A}_0+\mathcal{I}_{\lambda})/\mathcal{I}_{\lambda}$ both
dense in $\mathcal{A}/\mathcal{I}_{\lambda}$ and finite-dimensional.  Therefore
$$(\mathcal{A}_0+
\mathcal{I}_{\lambda})/\mathcal{I}_{\lambda}=\mathcal{A}/\mathcal{I}_{\lambda}.$$
Thus we have a homomorphism 
$S_{\lambda}=R_{\lambda}\circ \Gamma_{\lambda}\circ T_{\lambda}$ 
from $\mathbb{M}_2$ onto $\mathcal{A}/\mathcal{I}_{\lambda}$.
Since $\mathbb{M}_2$ has no non-trivial ideals,
the kernel of $S_{\lambda}$ is either $\mathbb{M}_2$ or $\{0\}$.  Since
$\mathcal{A}$ is a $C^*$-algebra, $\mathcal{I}_{\lambda}\neq \mathcal{A}$ 
(see \cite{BoSi}, p. 33),
and thus our homomorphism is injective; that is 
$\mathbb{M}_2\cong \mathcal{A}/\mathcal{I}_{\lambda}$.

Next consider (ii), with $\lambda={\bf p}$.  We repeat the above argument, but this
time, by Proposition~\ref{rangephi}, we may replace $\mathbb{M}_2$ on the left side
of (\ref{3hom}) by $\mathbb{M}_2^{diag}$.  Again, the above argument
yields a homomorphism $S_{{\bf p}}$ from $\mathbb{M}_2^{diag}$ onto
$\mathcal{A}/I_{{\bf p}}$.  However, unlike $\mathbb{M}_2$,
$\mathbb{M}_2^{diag}$ contains two non-trivial ideals, namely

\begin{equation}\label{possideals}
\left\{\left[\begin{array}{lr}
a&0\\0&0
\end{array}\right]:a\in\mathbb{C}\right\}\mbox{ and }
\left\{\left[\begin{array}{lr}
0&0\\0&b
\end{array}\right]:b\in\mathbb{C}\right\}.
\end{equation}
Again, $I_{{\bf p}}\neq \mathcal{A}$ and so $\mbox{ker }S_{{\bf p}}$ is either
$\{0\}$ or one of these two ideals.  If it is the first ideal 
in (\ref{possideals}),
then $S_{{\bf p}}$ induces an 
isomorphism of $\mathbb{C}$ and $\mathcal{A}/I_{{\bf p}}$
whose inverse has the form
$$[b]_{I_{{\bf p}}}\rightarrow w(\eta)$$
when $b$ is given by Equation~(\ref{defofa0}).  In particular, 
for $b=t_w$, we see
that $$\|[t_w]_{I_{{\bf p}}}\|=|w(\eta)|.$$
However, for $0<\lambda\leq s$, we know that
$$\|[t_w]_{\mathcal{I}_{\lambda}}\|=\left\|\left[\begin{array}{lr}
w(\zeta)&0\\0&w(\eta)\end{array}\right]\right\|_{\mathbb{M}_2}=
\max\{|w(\zeta)|,|w(\eta)|\}.$$
The map $\lambda\rightarrow\|[b]_{\mathcal{I}_{\lambda}}\|$ is known to be 
upper semi-continuous
on $\Lambda$ (see \cite{BoSi}, Theorem 1.34), which implies that 
for each $w$ in $C(\partial\D)$, 
$$
\max\{|w(\zeta)|,|w(\eta)|\} =
\limsup_{\lambda\downarrow 0}\|[t_w]_{\mathcal{I}_{\lambda}}\|
\leq  \|[t_w]_{I_{{\bf p}}}\|=|w(\eta)|.
$$
This is clearly impossible.  Thus $\mbox{ker }S_{{\bf p}}$ cannot be the first
ideal in (\ref{possideals}), or similarly, the second.  Therefore,
$S_{{\bf p}}$ has kernel $\{0\}$ and provides an isomorphism of $\mathbb{M}_2^{diag}$
and $\mathcal{A}/I_{{\bf p}}$, proving (ii).

Finally, for (iii), one can repeat the general argument from (i), with
$\lambda\in\partial\D\backslash\{\zeta,\eta\}$, replacing $\mathbb{M}_2$
in (\ref{3hom}) by $\{cI_{2\times 2}:d\in\mathbb{C}\}\cong \mathbb{C}$,
an algebra with no non-trivial ideals.  

One easily checks that the isomorphism $S_{\lambda}^{-1}$ 
from $\mathcal{A}/\mathcal{I}_{\lambda}$ into $\mathbb{M}_2$ is given for
$b$ in $\mathcal{A}_0$ by
$$ S_{\lambda}^{-1}:[b]_{\mathcal{I}_{\lambda}}\rightarrow \Phi_{\lambda}(b).$$
By Equation~(\ref{defphil}), $S_{\lambda}^{-1}$, and thus $S_{\lambda}$, 
are manifestly
$*-$maps. 
\end{proof}

\noindent {\bf Remark 2.} For future reference, we note that by the above
proof, the composition $S_{\lambda}$ of the three homomorphisms in (\ref{3hom}) 
is an
isomorphism, and thus the map 
$\Gamma_{\lambda}$
is an isomorphism of $\mathcal{A}_0/\mbox{ker }\Phi_{\lambda}$ and 
$\mathcal{A}_0/(\mathcal{A}_0
\cap \mathcal{I}_{\lambda})$.  In other words, 
$\mbox{ker }\Phi_{\lambda}=\mathcal{A}_0\cap \mathcal{I}_{\lambda}$.

By Proposition~\ref{rangephi} and Proposition~\ref{quotient}, we have
$*$-isomorphisms 
\begin{equation}
\mathcal{A}/\mathcal{I}_{\lambda}\cong\mathcal{A}_0/\mbox{ker }\Phi_{\lambda}\cong
\left\{\begin{array}{l}
\mathbb{M}_2 \mbox{ when }0<\lambda\leq s,\\
\mathbb{M}_2^{diag}\mbox{ when }\lambda = {\bf p},\\
\{cI_{2\times 2}:c\in\mathbb{C}\} \mbox{ when }
\lambda\in\partial\D\backslash\{\zeta,\eta\},
\end{array}
\right.
\end{equation}
the composition being $S_{\lambda}^{-1}$.  The objects on the 
right are $C^*-$algebras,
so that $S_{\lambda}^{-1}$ is isometric.
Thus, for $b\in \mathcal{A}_0$, 
\begin{equation}
\|[b]_{\mathcal{I}_{\lambda}}\|_{\mathcal{A}/\mathcal{I}_{\lambda}}=
\|\Phi_{\lambda}(b)\|=
\left\{\begin{array}{l}
\left\|\left[\begin{array}{lr}
w(\zeta)+g(\lambda)& h(\lambda)\\
k(\lambda)&w(\eta)+f(\lambda)\end{array}\right]\right\|\mbox{ if } 0<\lambda\leq s,
\\[1mm]
\\
\left\|\left[\begin{array}{lr}
w(\zeta)& 0\\
0&w(\eta)\end{array}\right]\right\|\mbox{ if } \lambda={\bf p},
\\[1mm]
\\
\left\|\left[\begin{array}{lr}
w(\lambda)& 0\\
0&w(\lambda)\end{array}\right]\right\|\mbox{ if } 
\lambda\in\partial\D\backslash\{\zeta,\eta\},
\end{array}
\right.
\end{equation}
the norm on the right being the operator norm in $\mathbb{M}_2$.

Now we write $B(\Lambda,\mathbb{M}_2)$ for the $C^*$-algebra
of all bounded functions $F$ from $\Lambda$ to $\mathbb{M}_2$, with
norm
$$\|F\|=\sup_{\lambda\in\Lambda}\|F(\lambda)\|_{\mathbb{M}_2}.$$
We can define a $*$-homomorphism $\Phi$ from $\mathcal{A}_0$ to $B
(\Lambda,\mathbb{M}_2)$ by letting $\Phi(b)$ be the function whose value at
$\lambda$ in $\Lambda$ is $\Phi_{\lambda}(b)$.  We write $\mathcal{D}$ for
the range of $\Phi$.  According to the above results and Douglas' theorem,
$\|b\|_{\mathcal{A}} = \sup_{\lambda\in\Lambda}\|\Phi_{\lambda}(b)\|$,
so that $\Phi$ is an isometric $*$-isomorphism of $\mathcal{A}_0$ onto 
the $*$-algebra
$\mathcal{D}$.  It is easy to verify that $\mathcal{D}$ consists of all 
$$F=\left[\begin{array}{lr}
f_{11}&f_{12}\\f_{21}&f_{22}
\end{array}\right]$$
in $B(\Lambda,\mathbb{M}_2)$
such that each $f_{ij}$ is continuous on $\{{\bf p}\} \cup (0,s)$
and $\partial\D\backslash\{\zeta,\eta\}$, $f_{12}$ and $f_{21}$
vanish at ${\bf p}$ and on $\partial\D\backslash\{\zeta,\eta\}$,
$f_{11}=f_{22}$ on $\partial\D\backslash\{\zeta,\eta\}$, while
$f_{11}({\bf p})=\lim_{\lambda\rightarrow \zeta}f_{11}(\lambda)$
and $f_{22}({\bf p})=\lim_{\lambda\rightarrow\eta}f_{22}(\lambda)$,
the
limits being taken as $\lambda\rightarrow\zeta$ or $\lambda\rightarrow
\eta$ through points in $\partial\D\backslash\{\zeta,\eta\}$.  One easily
checks that $\mathcal{D}$ is closed in $B(\Lambda,\mathbb{M}_2)$.
Since $\Phi$ is isometric, $\mathcal{A}_0$ is complete. Since $\mathcal{A}_0$
is dense in $\mathcal{A}$, we can close the circle to obtain the following
result.

\begin{prop}\label{a0=a}
The algebra $\mathcal{A}_0$ coincides with $\mathcal{A}$, and 
$\mbox{ker }\Phi=\mathcal{I}_{\lambda}$.
\end{prop}

Let us define two closed subspaces $\mathcal{M}$ and $\mathcal{N}$ in 
$\mathcal{A}$:
$$\mathcal{M}\equiv\{f(x^*x):f\in C_0([0,s])\}, \ 
\mathcal{N}\equiv\{g(xx^*):g\in C_0([0,s])\}.$$
We have already seen that $\mathcal{A}_0$ is an algebraic direct sum
of the closed subspaces $\{t_w:w\in C(\partial\D)\}$, $\mathcal{M},\mathcal{N},
u\mathcal{M}$ and $u^*\mathcal{N}$.  Since $\mathcal{A}_0=\mathcal{A}$,
a Banach space, we have the following corollary.

\begin{cor}\label{dsa}
As a Banach space, $\mathcal{A}=C^*(T_z,\cp)/\mathcal{K}$ has the
direct sum decomposition
$$\mathcal{A}=\{t_w:w\in C(\partial\D)\}\oplus\mathcal{M}\oplus\mathcal{N}\oplus
u\mathcal{M}\oplus u^*\mathcal{N}.$$
\end{cor}

In summary we have the following:
\begin{theorem}\label{5psi}
The map $\Phi$ is a $*$-isomorphism of $\mathcal{A}$ onto
$\mathcal{D}$.
\end{theorem}

\noindent {\bf Remark 3.} Given the form of the algebra ${\mathcal D}$, it is
not hard to show that every irreducible representation of
$C^*(T_z,\cp)/{\mathcal K}$ is unitarily equivalent either to
one of the two-dimensional representations $\Phi_{\lambda}$,
$\lambda$ in $(0,s]$, or to one of the scalar
representations $\ell_{\lambda}:b\rightarrow w(\lambda)$,
$\lambda$ in $\partial {\mathbb D}$, where
$b$ is given by Equation~(\ref{defofa0}).

\subsection{$C^*(T_z,\cp)$ REVISITED AND THE MAP $\Psi$}
Let $E$ and $F$ be the spectral projections of $\cp^*\cp$ 
and $\cp\cp^*$ respectively,
which are associated to their common essential spectrum $[0,s]$.  We have
$$\cp^*\cp=E\cp^*\cp E+(I-E)\cp^*\cp(I-E)$$
and 
$$\cp\cp^*=F\cp\cp^* F+(I-F)\cp\cp^*(I-F).$$
Notice that the second term on the right-hand side of each of 
these expressions
is a finite rank operator.  Thus if $f$ and $g$ are continuous on
$\sigma(\cp^*\cp)=\sigma(\cp\cp^*)$, then
\begin{equation}\label{18}
f(\cp^*\cp)=f(E\cp^*\cp E)+K_1, \ 
g(\cp\cp^*)=g(F\cp\cp^* F)+K_2
\end{equation}
for finite rank operators $K_1$ and $K_2$.  
Also note that the maps $f\rightarrow f(E\cp^*\cp E)$ and 
$g\rightarrow f(F\cp\cp^* F)$
are isometries from $C_0([0,s])$ onto closed subspaces $\mathfrak{M}$ and
$\mathfrak{N}$ in $C^*(T_z,\cp)$.  

\begin{theorem}\label{bsdirectsum}
As a Banach space, 
$C^*(T_z,\cp)$ is the direct sum of closed subspaces:
\begin{equation}\label{dsC}
C^*(T_z,\cp)=\{T_w:w\in C(\partial\D\}\oplus\mathfrak{M}\oplus\mathfrak{N}\oplus
U\mathfrak{M}\oplus U^*\mathfrak{N}\oplus\mathcal{K}.
\end{equation}
\end{theorem}
\begin{proof} 
Given $B\in C^*(T_z,\cp)$, the coset $b=[B]$ satisfies 
Equation~(\ref{defofa0})
for unique $w\in C(\partial\D)$ and $f,g,h$ and $k$ in $C_0([0,s])$.  
Since the coset map $B\rightarrow [B]$ is one-to-one when
restricted to each of the first five direct summands
(for example, $[Uh(\cp^*\cp)]=uh(x^*x)$), we see that
\begin{equation}\label{Bform}
B=T_w+f(E\cp^*\cp E)+g(F\cp\cp^* F)+Uh(E\cp^*\cp E)+U^*k(F\cp\cp^* F)+K
\end{equation}
for a unique compact operator $K$.  
\end{proof}

Now consider the map $\Psi:C^*(T_z,\cp)\rightarrow\mathcal{D}$ defined by
$\Psi(B)=\Phi([B])$.
Clearly
we have the following result.

\begin{theorem}\label{psi}
We have a short exact sequence of $C^*$-algebras,
$$0\rightarrow \mathcal{K}\stackrel{i}{\rightarrow} C^*(T_z,\cp)\stackrel{\Psi}
{\rightarrow}\mathcal{D}\rightarrow 0,$$
where $i$ is inclusion.
\end{theorem}

\subsection{THE DENSE SEMI-POLYNOMIAL SUBALGEBRA $\mathcal{P}$}
We write $\mathcal{P}$ for the dense non-commutative semi-polynomial $*-$ algebra
consisting of finite linear combinations of all $T_w$, $w$ in $C(\partial\D)$,
all words in $\cp$ and $\cp^*$, and all compact operators.
Every element of $\mathcal{P}$ has the form

\begin{equation}\label{P}
B=T_w+f(\cp^*\cp)+g(\cp\cp^*)+\cp p(\cp^*\cp)+\cp^*q(\cp\cp^*)+K,
\end{equation}
where $w$ is in $C(\partial\D)$, $f,g,p$ and $q$ are polynomials with 
$f(0)=0=g(0)$,
and
$K$ is compact.  Cutting $\cp^*\cp$ and $\cp\cp^*$ down by the spectral 
projections $E$ and $F$
respectively, we find 
\begin{eqnarray*}
B&=&T_w+f(E\cp^*\cp E)+g(F\cp\cp^* F)+UE\sqrt{\cp^*\cp}p(\cp^*\cp)E\\
&+ &U^*F\sqrt{\cp\cp^*}
q(\cp\cp^*)F+K',
\end{eqnarray*}
where we have absorbed each of the finite ranks arising from Equations~(\ref{18}) 
into the new compact operator $K'$.  By Theorem~\ref{bsdirectsum}, $B$ 
determines each of
the six summands here.  Since $f,g,p$ and $q$ are polynomials, and so are determined
by their restrictions to $[0,s]$, the decomposition of $B$ in Equation~(\ref{P}) is
unique.  Since $\cp^*\cp-sC_{\varphi\circ\sigma}$ and
$\cp\cp^*-sC_{\sigma\circ\varphi}$, are compact,  we see that Equation~(\ref{P}) becomes
$$B=T_w+A_1+A_2+A_3+A_4+K''$$
where $K''$ is compact, and $A_1,A_2,A_3,A_4$ are finite linear combinations 
of composition operators
whose associated self-maps of $\D$ are taken from
the respective lists
$(\varphi\circ\sigma)_{n_1}$, $(\sigma\circ\varphi)_{n_2}$,
$(\varphi\circ\sigma)_{n_3}\circ\varphi$, and $(\sigma\circ\varphi)_{n_4}\circ\sigma$,
for integers $n_1,n_2\geq 1$ and $n_3,n_4\geq 0$,
where $\tau_n$ denotes the $n^{th}$ iterate of the map $\tau$.
Since all of these self-maps are distinct, Corollary 5.17 
in \cite{KM} says the
corresponding composition operators are linearly independent modulo $\mathcal{K}$.
Thus the operator $B$ determines the coefficients in each of the sums 
$A_1,A_2,A_3,A_4$,
and $w$ and $K''$ as well.
We summarize these observations in the following theorem.

\begin{theorem}\label{descofP}
Every operator in $\mathcal{P}$ is a sum of a unique Toeplitz operator with
continuous symbol,
a unique compact operator
and a unique finite linear combination of composition operators with associated
disk maps 
taken from the set
$$\{(\varphi\circ\sigma)_{n_1},(\sigma\circ\varphi)_{n_2},(\varphi\circ\sigma)_{n_3}
\circ\varphi,(\sigma\circ\varphi)_{n_4}\circ\sigma\}
$$
where $ n_k\geq 1$ for $k=1,2$ and $
n_k\geq 0$ for $k=3,4$.
\end{theorem}

For an operator $B$ given by Equation~(\ref{P}), the matrix function
$\Psi(B)$ can properly be called the ``symbol of $B$".  In particular,
if $r$ is the function defined on $\Lambda$ by
$r(\lambda)=\sqrt{\lambda}$ for $0<\lambda\leq s$ and $r(\lambda)=0$ otherwise,
then
$$\Psi(\cp)=\left[\begin{array}{lr}
0&r\\0&0
\end{array}\right].
$$

\subsection{ESSENTIAL SPECTRA AND ESSENTIAL NORMS IN $C^*(T_z,\cp)$.}
\begin{theorem}\label{essspeccomp}
Let $B$ in $C^*(T_z,\cp)$ be given by Equation~(\ref{Bform}).  The essential
spectrum of $B$ is the union of $w(\partial\D)$ with the image of 
$$\frac{1}{2}[f(t)+w(\eta)+g(t)+w(\zeta)\pm 
\sqrt{(f(t)+w(\eta)-g(t)-w(\zeta))^2+4h(t)k(t)}]$$
as $t$ ranges over $[0,s]$.
\end{theorem}

\begin{proof}
By Theorem~\ref{5psi} or Theorem~\ref{psi}, the essential spectrum of $B$ is
$$\{z\in\mathbb{C}:\mbox{ det }(\Phi_{\lambda}([B])-zI_{2\times 2})=0
\mbox{ for some }\lambda\in\Lambda\}.$$
Evaluating this determinant via Equation~(\ref{defphil}) gives the desired result.
\end{proof}

We start with some examples of Theorem~\ref{essspeccomp} in which $w=0$.

\vspace{.2in}

\noindent {\bf Example 1.} The essential spectrum of the real part of $\cp$ is
the interval
$[-\sqrt{s}/2,\sqrt{s}/2]$, where $s=|\varphi'(\zeta)|^{-1}$.
This follows from using $f(t)=g(t)=0$ and $h(t)=k(t)=\sqrt{t}$ in
Theorem~\ref{essspeccomp} to see that
$$\sigma_e(\cp+\cp^*)=[-\sqrt{s},\sqrt{s}].$$

\vspace{.2in}

\noindent {\bf Example 2.} The essential spectrum of the self-commutator $[\cp^*,\cp]$
is
$[-s,s].$  This is obtained from Theorem~\ref{essspeccomp},
using $f(t)=t$, $g(t)=-t$, and $k(t)=h(t)=0$.  Similarly, the anti-commutator
$\cp^*\cp+\cp\cp^*$ has essential spectrum $[0,s]$.

\vspace{.2in}

\noindent {\bf Example 3.} Let $$B_1=C_{\varphi\circ\sigma}+C_{\sigma\circ\varphi}+\cp-
C_{\sigma},$$ so that $f(t)=t/s=g(t)$, $h(t)=\sqrt{t}$ and $k(t)=-\sqrt{t}/s$.
Then $\sigma_e(B_1)$ is the parabolic curve $y^2+iy$, $-1\leq y\leq 1$.

\vspace{.2in}

\noindent {\bf Example 4.} Let $$B_2=C_{\varphi\circ\sigma}-C_{\sigma\circ\varphi}+\frac{1}{2}\cp-
C_{\sigma},$$ so that $f(t)=t/s$, $g(t)=-t/s$, $h(t)=\sqrt{t}/2$ and 
$k(t)=-\sqrt{t}/s$.
Then $\sigma_e(B_2)$ is the union of two complex line segments, 
$[-\frac{1}{\sqrt{2}},\frac{1}
{\sqrt{2}}]$ and $[-\frac{i}{4},\frac{i}{4}]$.

\vspace{.2in}

\noindent {\bf Example 5.} Let $$B_3=2C_{\varphi\circ\sigma}+\cp-
C_{\sigma},$$ so that $f(t)=2t/s$, $g(t)=0$, $h(t)=\sqrt{t}$ and 
$k(t)=-\sqrt{t}/s$.
Here $\sigma_e(B_3)$ is the circle of radius $\frac{1}{2}$ centered 
at $z=\frac{1}{2}$.

\vspace{.2in}

Next we look at the effect of adding a Toeplitz operator.  Consider 
an operator $B=T_w+Y$
given by Equation~(\ref{Bform}), with
$$Y=f(E\cp\cp^*E)+g(F\cp\cp^*F)+Uh(E\cp^*\cp E)+U^*k(F\cp\cp^*F)+K.$$
According to Theorem~\ref{essspeccomp}, adding $Y$ to $T_w$ does not 
affect the part of
the essential spectrum coming from $\sigma_e(T_w)=w(\partial\D)$.  If
$w$ takes a common value $c$ at the points $\zeta$ and $\eta$, 
Theorem~\ref{essspeccomp}
also implies that 
$$\sigma_e(B)=\sigma_e(T_w)\cup\sigma_e(cI+Y).$$
In this case, the effect of adding $T_w$, on the part of the essential
spectrum coming from $Y$, is to merely translate it by $c$.  However if
$w(\zeta)\neq w(\eta)$, adding $T_w$ can non-trivially deform $Y's$ 
contribution
to $\sigma_e(B)$.  

\vspace{.2in}

\noindent {\bf Example 6.}  For $r\geq 0$, suppose $w$ in $C(\partial \D)$ satisfies
$$w(\eta)=r\frac{1+i}{\sqrt{2}}, w(\zeta)=-r\frac{1+i}{\sqrt{2}}.$$
Let $B=T_w+Y$ where $Y=\cp+\cp^*$.  Taking $f,g,h$, and $k$ as in Example 1,
we see from Theorem~\ref{essspeccomp} that 
$$\sigma_e(B)=w(\partial\D)\cup \{\pm\sqrt{t+r^2i}:0\leq t\leq s\}.$$
Thus when $r=0$ (so that $w(\zeta)=w(\eta)=0$), 
$$\sigma_e(B)=w(\partial\D)\cup[-\sqrt{s},\sqrt{s}]=\sigma_e(T_w)\cup\sigma_e(Y).$$
However, when  $r>0$, adding $T_w$ to $Y$ disconnects the essential spectrum
of the latter operator, deforming the two halves of 
$\sigma_e(Y),\  [0,\sqrt{s}]$ and
$[-\sqrt{s},0]$, into the curves $\{\sqrt{t+r^2i}:0\leq t\leq s\}$ and
$\{-\sqrt{t+r^2i}:0\leq t\leq s\}$, respectively.  The first of these curves
lies in the open first quadrant, is convex, and falls downhill to the right.
The second, of course, is its reflection through the origin.

\vspace{.2in}

Finally, we consider essential norms.  If $B$ in $C^*(T_z,\cp)$ is given by
Equation~(\ref{Bform}), we know that the essential norm $\|B\|_e$ is given by
$$\|B\|_e=\sup_{\lambda\in\Lambda}\|\Phi_{\lambda}([B])\|_{\mathbb{M}_2}.
$$

\vspace{.2in}

\noindent {\bf Example 7.} Let $B=T_z+\cp+\cp^*$.  Here we
have $w(e^{i\theta})=e^{i\theta}, f(t)=g(t)=0$ and $h(t)=k(t)=\sqrt{t}$.  If
$\lambda$ is in $\partial\D\backslash\{\zeta,\eta\}$ or $\lambda={\bf p}$, then
$\Phi_{\lambda}([B])$ is a diagonal unitary matrix.  For $0<\lambda\leq s$, 
$$\Phi_{\lambda}([B])= \left[\begin{array}{lr}
\zeta&\sqrt{\lambda}\\\sqrt{\lambda}&\eta
\end{array}\right].
$$ 
A well-known formula for the operator norm on $\mathbb{M}_2$ (see \cite{RR}, p.17)
gives
\begin{eqnarray*}
\|B\|^2_e&=&\sup_{0<\lambda\leq s}\left\|\left[\begin{array}{lr}
\zeta&\sqrt{\lambda}\\\sqrt{\lambda}&\eta
\end{array}\right]\right\|^2\\
&=&\sup_{0<\lambda\leq s}\left\{1+\lambda+
\sqrt{(1+\lambda)^2-|\zeta\eta-\lambda|^2}\right\}\\
&=&1+\frac{1}{|\varphi'(\zeta)|}
+\sqrt{\frac{2}{|\varphi'(\zeta)|}}\sqrt{1+\mbox{Re}(\zeta\eta)}.
\end{eqnarray*}

\end{document}